\documentstyle[12pt]{article}

\begin{document}
\newcommand{\ol }{\overline}
\newcommand{\ul }{\underline }
\newcommand{\ra }{\rightarrow }
\newcommand{\lra }{\longrightarrow }
\newcommand{\ga }{\gamma }
\newcommand{\st }{\stackrel }
\newcommand{\scr }{\scriptsize }
\title{\Large\textbf{ On Polynilpotent Multipliers of Free Nilpotent
Groups}\footnote{This research was in part supported by a grant from
IPM (No. 84200038).}}
\author{\textbf{Mohsen Parvizi\footnote{Email address:
parvizi@dubs.ac.ir}\ \ and
  Behrooz Mashayekhy
\footnote{Correspondence: mashaf@math.um.ac.ir}} \\
Department of Mathematics, Ferdowsi University of Mashhad,\\
P. O. Box 1159-91775, Mashhad, Iran\\
 and\\
 Institute for Studies in Theoretical Physics and Mathematics (IPM)  }
\date{ }

\maketitle
\begin{abstract}
In this paper, we present an explicit structure for the Baer
invariant of a free nilpotent group ( the $n$-th nilpotent product
of the infinite cyclic group, $\textbf{ Z}\st{n}* \textbf{
Z}\st{n}*\ldots \st{n}*\textbf{ Z}$) with respect to the variety
of polynilpotent groups of class row $(c,1)$, ${\cal N}_{c,1}$,
for all $c > 2n-2$. In particular, an explicit structure of the
Baer invariant of a free abelian group with respect to the
variety of metabelian groups will be presented.
\end{abstract}

\hspace{-0.6cm}\textit{2000 Mathematics Subject Classification }:20E34, 20E10, 20F18\\
\hspace{-0.6cm}\textit{Keywords}: Baer invariant, free nilpotent
group, polynilpotent variety, nilpotent product\\
\newpage
\begin{center}
\hspace{-0.65cm}\textbf{1.Introduction and Preliminaries}\\
\end{center}
Historically, there have been several papers from the beginning of
twentieth century, trying to find some structure for the well-known
notion of the Schur multiplier and its varietal generalization, the
Baer invariant, of some famous products of groups, such as direct
products, free products, nilpotent products and regular products.\\

I. Schur [13] in 1907 and J. Wiegold [14] in 1971 obtained the
structure of the Schur multiplier of the direct product of two
finite groups as follows:
\begin{center}
$M(A\times B)\cong M(A)\oplus M(B)\oplus \frac{[A,B]}{[A,B,A*B]}$,
where $\frac{[A,B]}{[A,B,A*B]}\cong A_{ab}\otimes B_{ab}$.
\end{center}

In 1979, M. R. R. Moghaddam [11] and in 1998, G. Ellis [2], tried to
generalize the above result to obtain the structure of the
$c$-nilpotent multiplier of the direct product of two groups, ${\cal
N}_cM(A\times B)$. Also in 1997 the second author in a joint paper
[7] presented an explicit formula for the $c$-nilpotent
multiplier of a finite abelian group.\\

In 1972 W. Haebich [4] presented a formula for the Schur
multiplier of a regular product of a family of groups. It is known
that the regular product is a generalization of the nilpotent
product and the last one is a generalization of the direct
product, so Haebich's result is an interesting generalization of
the Schur-Wiegold result. Also, M. R. R. Moghaddam [12], in 1979
gave a formula similar to Haebich's formula for the Schur
multiplier of a nilpotent product. Moreover, in 1992, N. D. Gupta
and M. R. R. Moghaddam [3] tried to present an explicit formula
for the $c$-nilpotent multiplier of the $n$-th nilpotent product
$\textbf{Z}_2\st{n}*\textbf{Z}_2$ (see [8, Defn. 2.6] for the
definition of the $n$-th nilpotent
product).\\

In 2001, the second author [8] found a structure similar to
Haebich's type for the $c$-nilpotent multiplier of a nilpotent
product of a family of cyclic groups. The $c$-nilpotent multiplier
of a free product of some cyclic groups was studied by the second
author [9] in 2002.\\
Finally the authors [10] concentrated on the Baer invariant with
respect to the variety of polynilpotent groups, for the first
time. We succeeded to present an explicit formula for the
polynilpotent multipliers of finitely generated abelian groups.\\
Now in this paper we intend to obtain an explicit formula for some
polynilpotent multipliers of an $n$-th nilpotent product of some
infinite cyclic groups,
\begin{center}
${\cal N}_{c,1}M(\textbf{ Z}\st{n}* \textbf{ Z}\st{n}*\ldots
\st{n}*\textbf{ Z})$, for all $c> 2n-2$.
\end{center}
(Note that $\textbf{ Z}\st{n}* \textbf{ Z}\st{n}*\ldots
\st{n}*\textbf{ Z}$  has the presentation $F/\ga_{n+1}(F)$, where $F$ is a free
group and hence it is a free nilpotent group of class $n$).\\

In particular, the structure of the metabelian multiplier of
direct product of some infinite cyclic group, ${\cal
S}_2M(\textbf{ Z}\oplus \textbf{ Z}\oplus \ldots \oplus\textbf{
Z})$ is completely known.\\

In the following we present some preliminaries which are used in
our method.\\

\hspace{-0.65cm}\textbf{Definition 1.1.} Let $G$ be any group with
a free presentation $G\cong F/R$, where $F$ is a free group. Then,
after R. Baer [1], the \textit{Baer invariant} of $G$ with respect
to a variety of groups ${\cal V}$, denoted by ${\cal V}M(G)$, is
defined to be
$${\cal V}M(G)=\frac{R\cap V(F)}{[RV^*F]}\ ,$$
where $V$ is the set of words of the variety ${\cal V}$, $V(F)$ is
the verbal subgroup of $F$ with respect to ${\cal V}$ and
$$[RV^*F]=<v(f_1,\ldots ,f_{i-1},f_ir,f_{i+1},\ldots,f_n)v(f_1,\ldots,f_i,
\ldots,f_n)^{-1}\mid $$
$$r\in R, 1\leq i\leq n, v\in V ,f_i\in F, n\in {\bf N}>.$$

In special case, if ${\cal V}$ is the variety of abelian groups,
${\cal A}$ , then the Baer invariant of $G$ will be the well-known
notion the \textit{Schur multiplier}
$$\frac{R\cap F'}{[R,F]}.$$

If ${\cal V}$ is the variety of nilpotent groups of class at most
$c\geq1$, ${\cal N}_c$, then the Baer invariant of $G$ with
respect to ${\cal N}_c$ which is called the \textit{$c$-nilpotent
multiplier} of $G$, will be
$${\cal N}_cM(G)=\frac{R\cap \gamma_{c+1}(F)}{[R,\ _cF]},$$
where $\gamma_{c+1}(F)$ is the $(c+1)$-st term of the lower
central series of $F$ and $[R,\ _1F]=[R,F], [R,\ _cF]=[[R,\
_{c-1}F],F]$, inductively.

 In a very more general case, let $\cal V$ be the variety of
 polynilpotent groups of class row $(c_1,\ldots ,c_t)$, ${\cal
 N}_{c_1,\ldots ,c_t}$, then the Baer invariant of a group $G$
 with respect to this variety, which we call it a \textit{polynilpotent multiplier}, is as
 follows:
 $${\cal N}_{c_1,\ldots ,c_t}M(G)\cong \frac{R\cap
 \ga_{c_t+1}\circ\ldots\circ\ga_{c_1+1}(F)}{[R{\cal N}_{c_1,\ldots ,c_t}^*F]}$$
 $$\cong
 \frac{R\cap \ga_{c_t+1}\circ\ldots\circ\ga_{c_1+1}(F)}{[R,\ _{c_1}F,\ _{c_2}\ga_{c_1+1}(F),
 \ldots , \ _{c_t}\ga_{c_{t-1}+1}\circ \ldots \circ
 \ga_{c_1+1}(F)]},\ \ (\star)$$
 where
 $\ga_{c_t+1}\circ\ldots\circ\ga_{c_1+1}(F)=\ga_{c_t+1}(\ga_{c_{t-1}+1}
 (\ldots(\ga_{c_1+1}(F))\ldots
 ))$ is a term of the iterated lower central series of $F$ which is the verbal subgroup
 of $F$ with respect to ${\cal
 N}_{c_1,\ldots ,c_t}$.
See [6, Corollary 6.14] for the equality
 $$[R{\cal N}_{c_1,\ldots ,c_t}^*F]=[R,\ _{c_1}F,
 \ _{c_2}\ga_{c_1+1}(F),\ldots ,
\ _{c_t}\ga_{c_{t-1}+1}\circ \ldots \circ \ga_{c_1+1}(F)].$$

\hspace{-0.7cm}\textbf{Definition 1.2.}\ \textit{Basic
commutators} are defined in the usual way. If $X$ is a fully
ordered independent subset of a free group, the basic
commutators on $X$ are defined inductively over their weight as follows:\\
$(i)$ All the members of $X$ are basic commutators of weight one on $X$;\\
$(ii)$ Assuming that $n>1$ and that the basic commutators of
weight less than $n$ on $X$ have been defined and ordered;\\
$(iii)$ A commutator $[a,b]$ is a basic commutator of weight $n$
on $X$ if  $wt(a)+wt(b)=n,\ a<b$, and if $b=[b_1,b_2]$, then
$b_2\leq a$. The ordering of basic commutators is then extended to
include those of weight $n$ in any way such that those of weight
less than $n$ precede those of weight $n$. The natural way to
define the order on basic commutators of the same weight is
lexicographically, $[b_1,a_1]<[b_2,a_2]$ if $b_1<b_2$ or if
$b_1=b_2$ and $a_1<a_2$.

The next two theorems are vital in our investigation.\\

\hspace{-0.65cm}\textbf{Theorem 1.3} (P.Hall [5]).\textit{ Let
$F=<x_1,x_2,\ldots ,x_d>$ be a free group, then
$$ \frac {\ga_n(F)}{\ga_{n+i}(F)} \ \ , \ \ \ \  1\leq i\leq n$$
is the free abelian group freely generated by the basic
commutators of weights
$n,n+1,\ldots ,n+i-1$ on the letters $\{x_1,\ldots ,x_d\}.$}\\

\hspace{-0.65cm}\textbf{Theorem 1.4} (Witt Formula [5]).\textit{
The number of basic commutators of weight $n$ on $d$ generators is
given by the following formula:
$$ \chi_n(d)=\frac {1}{n} \sum_{m|n}^{} \mu (m)d^{n/m}\ \ ,$$
where $\mu (m)$ is the {\it Mobious function}, which is defined to
be
   \[ \mu (m)=\left \{ \begin{array}{ll}
      1 & {\rm if}\ m=1, \\ 0 & {\rm if}\ m=p_1^{\alpha_1}\ldots
p_k^{\alpha_k}\ \ \exists \alpha_i>1, \\ (-1)^s & {\rm if}\
m=p_1\ldots p_s,
\end{array} \right.  \]  \\
where the $p_i$, $1\leq i\leq k$, are the distinct primes dividing
$m$.}
\begin{center}
\hspace{-0.65cm}\textbf{2  The Main Results}\\
\end{center}
Let $G\cong \textbf{ Z}\st{n}* \textbf{ Z}\st{n}*\ldots
\st{n}*\textbf{ Z}$ ($m$-copies of $\textbf{ Z}$) be the $n$-th
nilpotent product of $m$ copies of the infinite cyclic group
$\textbf{ Z}$. It is known that $G$ is the free $n$-th nilpotent
group of rank $m$ or equivalently $G$ is the free nilpotent group of
class $n$, and so has the following free presentation
$$1\lra \ga_{n+1}(F)\lra F\lra G\lra 1,$$
where $F$ is a free group on a set $X=\{x_1,x_2,\ldots,x_m\}$.\\
Clearly the Baer invariant of $G$ with respect to the variety of
nilpotent groups of class at most $c\geq 1$, ${\cal N}_c$, is
defined as follows
$${\cal N}_cM(G)\cong \frac{\ga_{n+1}(F)\cap \ga_{c+1}(F)}{[\ga_{n+1}(F),\ _c\ F]}=\frac{
\ga_{n+1}(F)\cap \ga_{c+1}(F)}{\ga_{n+c+1}(F)}.$$
In the following theorem, the structure of ${\cal N}_cM(G)$ will
be presented.\\

\hspace{-0.7cm}\textbf{ Theorem 2.1.} \textit{Let $G$ be a free $n$-th nilpotent group of rank
$m$. Then\\
$(i)\ {\cal N}_cM(G)$ is the free abelian group of rank
$\st{c+n}\sum_{i=c+1} \chi_i(m)$, for all
$c\geq n$.\\
$(ii) {\cal N}_cM(G)$ is the free abelian group of rank $\st{c+n}\sum_{i=n+1} \chi_i(m)$, for
 all $c\leq n$.}\\

\hspace{-0.7cm}\textit{ Proof.} $(i)$ Let $c\geq n$. Then
$${\cal N}_cM(G)\cong \frac{\ga_{c+1}(F)}{\ga_{c+n+1}(F)}.$$
Now by P. Hall's Theorem (1.3) $\ga_{c+1}(F)/\ga_{c+n+1}(F)$ is a
free abelian group freely generated by all basic commutators of
weights $c+1,\ldots,c+n$ on $X$. So its rank is
$\st{c+n}\sum_{i=c+1}
\chi_i(m)$.\\
$(ii)$ Since $c\leq n$, ${\cal N}_cM(G)\cong
\ga_{n+1}(F)/\ga_{c+n+1}(F).$ Hence the result holds
similar to $i$. $\Box$\\

Now, we try to obtain the structure of some polynilpotent
multipliers of $G$ of the form
$${\cal N}_{c,1}M(G),$$
where $c> 2n-2$. Using $(\star)$ we have
$${\cal N}_{c,1}M(G)\cong \frac{\ga_{n+1}(F)\cap \ga_2(\ga_{c+1}(F))}{[\ga_{n+1}(F),\ _c\ F,
\ga_{c+1}(F)]}=\frac{\ga_2(\ga_{c+1}(F))}{[\ga_{n+c+1}(F),\ga_{c+1}(F)]}.$$
In order to find the structure of ${\cal N}_{c,1}M(G)$, we need
the following notation and lemmas. Using Definition and Notation 1.2, put \\
$Y$= The set of all basic commutators on X of weights $c+1,\ldots
,c+n$\\
$Z$= The set of all basic commutators on Y of weight 2.\\

\hspace{-0.7cm}\textbf{ Lemma 2.2.}\textit{ With the above
notation, if $c\geq n-1$, then every element of $Z$
is a basic commutator on $X$}.\\

\hspace{-0.7cm}\textit{ Proof.} Every element of $Z$ has the form
$[b,a]$, where $b$ and $a$ belong to $Y$ and $b>a$ $i.e.$ $b$ and
$a$ are basic commutators of weights $c+i$ and $c+j$ on $X$ in
which $1\leq j\leq i\leq n$,
respectively.\\
Now, let $b=[b_1,b_2]$. In order to show that $[b,a]$ is a basic
commutator on $X$, it is enough to show that $b_2\leq a$. Since
$b=[b_1,b_2]$ is a basic commutator on $X$, so $b_1>b_2 $ and
hence $wt(b_2)\leq \frac{1}{2}wt(b)$, since $b\in Y$, so
$wt(b)\leq c+n$. Now, if $c\geq n-1$, then
$\frac{1}{2}(c+n)<c+1$. Thus, we have
$$wt(b_2)\leq \frac{1}{2}wt(b)\leq \frac{1}{2}(c+n)<c+1\leq c+i=wt(a).$$
Therefore $b_2<a$ and hence the result holds. $\Box$\\

\hspace{-0.7cm}\textbf{ Lemma 2.3.} \textit{With the above
notation and assumption, if $c\geq n-1$, then we have}
$$\ga_2(\ga_{c+1}(F))\equiv <Z> \ \ \ mod \ \ [\ga_{c+n+1}(F),\ga_{c+1}(F)]. $$

\hspace{-0.7cm}\textit{ Proof.} Let $[\alpha,\beta]$ be a
generator of $\ga_2(\ga_{c+1}(F))$, so ${\alpha},{\beta}\in
\ga_{c+1}(F)$. By P. Hall's Theorem (1.3) we can put
$\alpha=\alpha_1\alpha_2\ldots\alpha_t\eta$ \ and \
$\beta=\beta_1\beta_2\ldots\beta_s\mu$, where
$\alpha_1,\ldots,\alpha_t$ and $\beta_1,\ldots,\beta_s$  are basic
commutators of weights $c+1,\ldots,c+n$ on $X$ and
$\eta\ , \ \mu \in\ga_{c+n+1}(F)$.\\
By using commutator calculus, it is easy to see that
$$[\alpha,\beta]=\Pi_{i,j} [\alpha_i,\beta_j]^{f_{ij}}[\alpha_i,\mu]^{g_i}[\eta,\beta_j]^{h_j},$$
where $f_{ij},g_i,h_j\in\ga_{c+1}(F).$ Note that, since $c\geq
n-1$, $wt([\alpha_i,\beta_j])\geq 2c+2\geq c+n+1$. Now, it is
easy to see that
\begin{center}
$[\alpha_i,\mu],[\eta,\beta_j],[\alpha_i,\beta_j,f_{ij}]\in
[\ga_{c+n+1}(F),\ga_{c+1}(F)]$, for all $1\neq f_{ij}\in
\ga_{c+1}(F)$.
\end{center}
Therefore we have
$$[\alpha,\beta]=\Pi_{i,j}[\alpha_i,\beta_j] \ \ (mod [\ga_{c+n+1}(F),\ga_{c+1}(F)]).$$
It is easy to see that $\Pi_{i,j}[\alpha_i,\beta_j]\in <Z>$. Hence
the result holds.$\Box$\\

Now, we define the following useful set
\begin{center}
$W=\{ \ [b,a] \ | \ b$ and $a$ are basic commutators on $X$ such
that $b>a$, $wt(b)\geq c+n+1$ , $wt(a)\geq c+1$ , $wt(b)+wt(a)\leq
2c+2n+1 \}$.
\end{center}

\hspace{-0.7cm}\textbf{ Lemma 2.4.} \textit{With the above
notation, if $c> 2n-2$,
then every element of $W$ is a basic commutator on $X$.}\\

\hspace{-0.7cm}\textit{ Proof.} Let $[b,a]$ be an element of $W$.
By definition of $W$, $b$ and $a$ are basic commutators on $X$ and
$b>a$. Now, let $b=[b_1,b_2]$, where $b_1$ and $b_2$ are basic
commutators on $X$. It is enough to show that $b_2\leq a$. Since
$b=[b_1,b_2]$ is a basic commutator, so $b_1>b_2$ and hence
$wt(b_2)\leq \frac{1}{2}wt(b)$. By definition of $W$, $c+1\leq
wt(a), \ c+n+1\leq wt(b)\leq 2c+2n+1-wt(a)\leq
2c+2n+1-(c+1)=c+2n$. Since $c> 2n-2$ we have $\frac{1}{2}(c+2n)<
c+1$. Thus $wt(b_2)\leq \frac{1}{2}wt(b)\leq \frac{1}{2}(c+2n)<
c+1\leq
wt(a).$ Hence the result holds.$\Box$\\

\hspace{-0.7cm}\textbf{ Lemma 2.5.} \textit{With the above
notation and assumption, if $c\geq n-1$ then we have}
$$[\ga_{c+n+1}(F),\ga_{c+1}(F)]\equiv <W> \ \ mod \ \ \ga_{2c+2n+2}(F).$$

\hspace{-0.7cm}\textit{ Proof.} Let $[\alpha,\beta]$ be a
generator of $[\ga_{c+n+1}(F),\ga{c+1}(F)]$, so $\alpha\in
\ga_{c+n+1}(F)$ and $\beta\in \ga_{c+1}(F)$. By P. Hall's Theorem
(1.3) and considering two free abelian groups
$\ga_{c+n+1}(F)/\ga_{c+2n+1}(F)$
and $\ga_{c+1}(F)/\ga_{c+n+1}(F)$ we can write\\
$\alpha=\alpha_1\alpha_2\ldots\alpha_t\eta$ and
$\beta=\beta_1\beta_2\ldots\beta_s\mu$, where
$\alpha_1,\ldots,\alpha_t$ are basic commutators of weights
$c+n+1,\ldots,c+2n$ on $X$ and $\eta\in\ga_{c+2n+1}(F)$, and
$\beta_1,\ldots,\beta_s$ are basic commutators of weights
$c+1,\ldots,c+n$ on $X$ and $\mu\in\ga_{c+n+1}(F)$.\\
By using commutator calculus, it is easy to see that
$$[\alpha,\beta]=\Pi_{i,j}[\alpha_i,\beta_j]^{f_{ij}}[\alpha_i,\mu]^{g_i}[\eta,\beta_j]^{h_j},$$
where $f_{ij},\ g_i,\ h_j\in \ga_{c+1}(F).$ Now we have
$$wt(\alpha_i)+wt(\mu)\geq(c+n+1)+(c+n+1)\geq 2c+2n+2 $$
$$wt(\eta)+wt(\beta_j)\geq (c+2n+1)+(c+1)\geq 2c+2n+2 $$
\begin{center}
$wt(\alpha_i)+wt(\beta_j)+wt(f_{ij})\geq
(c+n+1)+(c+1)+(c+1)=3c+n+3\geq 2c+2n+2$  for all $1\neq f_{ij}\in
\ga_{c+1}(F)$, since $c\geq n-1$.
\end{center}
Therefore $$[\alpha,\beta]\equiv \Pi_{i,j}[\alpha_i,\beta_j]\ \
mod \ \ \ga_{2c+2n+2}(F) $$
$$ \equiv \Pi_{wt(\alpha_i)+wt(\beta_j)\leq 2c+2n+1}[\alpha_i,\beta_j]\ \ mod \ \
\ga_{2c+2n+2}(F)\in W. $$ Hence the result holds.$\Box$\\

Now, we are in a position to state and proof the main result of
the paper.\\

\hspace{-0.7cm}\textbf{ Theorem 2.6} \textit{Let $\textbf{
Z}\st{n}* \textbf{ Z}\st{n}*\ldots \st{n}*\textbf{ Z} \cong
F/\ga_{n+1}(F)$ be the free nilpotent group of class $n$. With
the above notation and assumption, if $c> 2n-2$, then ${\cal
N}_{c,1}M(\textbf{ Z}\st{n}* \textbf{ Z}\st{n}*\ldots
\st{n}*\textbf{ Z})$ is a free abelian group with the following
basis
$$B=\{ \ z[\ga_{c+n+1}(F),\ga_{c+1}(F)]\ \ | \ \ z\in Z \ \}$$
(Note that if $ c > 2n-2 $, since $ c \geq 1 $, then $ c \geq n-1 $)}.\\

\hspace{-0.7cm}\textit{ Proof.} Clearly ${\cal
N}_{c,1}M(\underbrace{\textbf{ Z}\st{n}* \textbf{ Z}\st{n}*\ldots
\st{n}*\textbf{
Z}}_{m-copies})=\ga_2(\ga_{c+1}(F))/[\ga_{c+n+1}(F),\ga_{c+1}(F)]$
is an abelian group which is generated by $B$, using Lemma 2.3. So
it is enough to show that elements of $B$ are linearly independent.
Suppose $\bar{z_1},\ldots,\bar{z_k}\in B$ and
$\sum_{i=1}^{k}\alpha_i\bar{z_i}=0$, where $\alpha_i$'s are
integers. Clearly $\sum_{i=1}^{k}\alpha_iz_i\in
[\ga_{c+n+1}(F),\ga_{c+1}(F)].\ \ \ (\star\star)$\\
Now consider the abelian group $\ga_{2c+2}(F)/\ga_{2c+2n+2}(F)$,
which is free abelian with the basis of all basic commutators on
$X$ of weights $2c+2,\ldots,2c+2n+1$, by P. Hall's Theorem (1.3).
By Lemma 2.5 we have
$$\frac{\ga_{2c+2n+2}(F)[\ga_{c+n+1}(F),\ga_{c+1}(F)]}{\ga_{2c+2n+2}(F)}=< \ w\ga_{2c+2n+2}(F)\
\ | \ \ w\in W \ >.$$ By $(\star\star)$, there exists
$w_1,\ldots,w_t\in W$, and integer numbers
$\beta_1,\ldots,\beta_t$ such that
$$\sum_{i=1}^{k}\alpha_i(z_i\ga_{2c+2n+2}(F))=\sum_{i=1}^{t}\beta_i(w_i\ga_{2c+2n+2}(F)).$$
Therefore
$$\sum_{i=1}^{k}\alpha_i(z_i\ga_{2c+2n+2}(F))+\sum_{i=1}^{t}(-\beta_i)(w_i\ga_{2c+2n+2}(F))=0 .\
 \ \ (\star\star\star)$$
By Lemmas 2.2 and 2.4 every element of $Z$ and $W$ is a basic
commutator on $X$ of weights $2c+2,\ldots,2c+2n+1$. Also by
considering the form of elements of $Z$ and $W$, it is easy to see
that $Z\cap W=\phi$. Now, by $(\star\star\star)$ and considering
the basis of the free abelian group
$\ga_{2c+2}(F)/\ga_{2c+2n+2}(F)$ we have $\alpha_i=0$ for all
$1\leq i \leq k$ and $\beta_i=0$ for all
$1\leq i \leq t$. Hence the result holds.$\Box$\\

\hspace{-0.7cm}\textbf{ Corollary 2.7.} \textit{If $c> 2n-2$, then\\
${\cal N}_{c,1}M(\underbrace{\textbf{ Z}\st{n}* \textbf{
Z}\st{n}*\ldots \st{n}*\textbf{ Z}}_{m-copies})\cong \textbf{
Z}\oplus \ldots \oplus\textbf{ Z} \ \ \
(\chi_2(\chi_{c+1}(m)+\ldots+\chi_{c+n}(m))-$copies), where
$\chi_j(i)$ is the number of all basic commutators of weight
$j$ on $i$ letters.\\
In particular \\
${\cal S}_2M(\underbrace{\textbf{ Z}\oplus \textbf{ Z}\oplus\ldots
\oplus\textbf{ Z}}_{m-copies})\cong \textbf{
Z}\oplus\ldots \oplus \textbf{ Z} \ \ \ (\chi_2(\chi_2(m))$-copies),\\
where ${\cal S}_2$ is the variety of all metabelian groups.
(solvable groups of length at most 2). Note that the authors
presented a similar structure for ${\cal S}_2M({\bf Z}\oplus {\bf
Z}\oplus\ldots\oplus {\bf Z})$ in a different method. (See [10])}.

\end{document}